\documentclass[smallextended]{svjour3}
\usepackage{latexsym,amssymb,amsmath}
\smartqed
\usepackage{color}
\usepackage{hyperref}
\RequirePackage{fix-cm}
\usepackage[section]{placeins}

\def\F{Fr\'{e}chet}
\def\O{L\'{o}pez}
\def\gph{\mbox{\rm gph}\,}
\def\epi{{\rm epi}\,}

\def\dom{\mbox{\rm dom}\,}

\begin{document}
	
	\renewcommand{\theequation}{\thesection.\arabic{equation}}
	\normalsize
	
	\setcounter{equation}{0}
	\title{Sensitivity Analysis  in Parametric  Convex
 Vector Optimization}
	
	\author{Duong Thi Viet An       \and
		Le Thanh Tung
	}

	\institute{Duong Thi Viet An\ \at
		Department of Mathematics and Informatics, Thai Nguyen University of Sciences, Thai Nguyen city 250000, Vietnam\\
		\email{andtv@tnus.edu.vn}           	\and
		Le Thanh Tung
		 \at
	Department of Mathematics, College of Natural Sciences, Can Tho University, Can Tho 900000, Vietnam\\
		\email{lttung@ctu.edu.vn}
	}
	\date{Received: date / Accepted: date}
	
	\maketitle
	
	\begin{abstract}
		In this paper, sensitivity  analysis of the efficient sets in parametric convex vector optimization is considered. Namely, the perturbation, weak perturbation, and proper perturbation maps are defined as set-valued maps. We establish the formulas for computing the Fr\'{e}chet coderivative of the profile of the above three kinds of perturbation maps. Because of the convexity assumptions, the conditions set are fairly simple if compared to those in the general case. In addition, our conditions are stated directly on the data of the problem. It is worth emphasizing that our approach is based on convex analysis tools which are different from those in the general case.
	\end{abstract}
	

	
	\newpage
\section{Introduction}
\markboth{\centerline{\it Sensitivity Analysis  in Parametric  Convex
		Vector Optimization}}{\centerline{\it D.T.V.~An and L.T.~Tung}} \setcounter{equation}{0}
\textit{Stability and sensitivity analysis} of parametric optimization problems are not only theoretically interesting but also practically important in optimization and variational analysis. Usually, by \textit{stability analysis} we mean the qualitative analysis, that is, the study of various continuity properties of the perturbation (or marginal) function (or map) of a family of parameterized optimization problems. Meanwhile, by \textit{sensitivity analysis} (or differential stability) we mean the quantitative analysis, that is, the study of derivatives (coderivative) of the perturbation function. We refer the interested reader to the books by Fiacco~\cite{Fiacco_1983}, by Malanowski~\cite{Malanowski_1987} and by Bonnans and Shapiro~\cite{Bonnans_Shapiro_2000} for a systematic view of these large topics.

 For a vector optimization problem that depends on a parameter vector, the  stability and sensitivity analysis of perturbation, weak perturbation, and proper perturbation maps are dealt with. Each of the perturbation maps is defined as a set-valued map that associates to each parameter value the set of all minimal, weakly minimal, and properly minimal points of the perturbed feasible set in the objective space with respect to a fixed ordering cone. 
Various stability and sensitivity analysis results in this direction can be found in the books~\cite[Chapter 4]{Sawaragi_Nakayama_Tanino_85},~\cite[Chapter 4]{Luc_book_89},~\cite[Chapter~13]{Khan_Tammer_Zanilescu_2015} and  in the old and new papers~\cite{C11,Chuong_Huy_Yao,CY09,CY10,CY13,HMY08,Kuk_Tanino_Tanaka_96,LH07,S93,Tanino_1988a,Tanino_1988b,Tanino_1990,Tung2016,Tung2020,Tung_2021}.

 In this paper,  sensitivity  analysis of the perturbation maps in convex vector optimization is considered. More precisely, we study the behavior of the perturbation, weak perturbation, and proper perturbation maps. As mentioned above, these maps are set-valued maps. Thus, in order to obtain the differential property of these maps, one may need to evaluate \textit{generalized derivatives} for set-valued maps.
 
 In the literature on differentiability properties of the perturbation maps, the papers of Tanino~\cite{Tanino_1988a,Tanino_1988b,Tanino_1990} are among the first results in this topic. Namely, the behavior of the perturbation map is analyzed quantitatively by using the concept of \textit{contingent derivatives} (\textit{codifferential}) for set-valued maps introduced by Aubin and Ekeland in~\cite[Chapter 4]{Aubin_Ekeland_1984}. Later, by using the concepts of contingent derivatives and \textit{Dini derivatives}, the authors in~\cite{Kuk_Tanino_Tanaka_96} investigate the quantitative information about the behavior of the above three kinds of perturbation maps. Another approach applied by Taa in~\cite{Taa_2003,Taa_2005} and An and Guti\'{e}rrez in~\cite{An_Cesar_2021} is to use the (weak) \textit{subdifferential} for multifunctions.
 
 It is worth emphasizing that the contingent derivatives mentioned above for set-valued maps are generated by \textit{tangent cones} to their graphs and it is known as the \textit{approximation in primal spaces}. The other line of the graphical approach to generalized differentiation was developed by Mordukhovich in~\cite{Mordukhovich_1980}, the \textit{coderivative} notion for general set-valued mappings via the normal cones to their graphs (the \textit{approximation in dual spaces}).  This is conceptually different from tangentially generated derivatives in the line of Aubin and Ekeland due to the absence of duality between tangent and normal cones in general nonconvex settings. Obviously, for smooth and convex-graph maps the two approaches are equivalent. By using these concepts of coderivative for multifunctions, the researchers have obtained some sensitivity analysis results 
for a vector optimization problem, see~\cite{C11,CY09,CY10,CY13,HMY08} and the references therein.

The primary objective of this paper is to obtain formulas for computing the \textit{Fr\'{e}chet coderivative of the  profile} of perturbation, weak perturbation, and proper perturbation maps in parametric  convex vector optimization. The first main result is presented in Theorem~\ref{Thm_frontier_map}. Namely, we obtain the exact formulas for computing the Fr\'{e}chet coderivative of the profile of the perturbation map. This proof is divided into two parts. First, we compute the Fr\'{e}chet coderivative of the profile of the objective map as in Theorem~\ref{Thm_objective_map}. Because of the convexity assumptions, the profile of the objective map is convex, it follows that this is the coderivative of the convex set-valued map. The important point to note here is the perturbation maps are nonconvex in general.
In the second part, by using the \textit{domination property}, we get the corresponding Fr\'{e}chet coderivative of the profile of the perturbation map. By some suitable changes,  the results in this paper still hold for the weak perturbation and proper perturbation maps (see Theorems~\ref{Weak_theorem} and~\ref{Proper_theorem}).

The Fr\'{e}chet coderivative of the  profile of the perturbation map is known as the \textit{Fr\'{e}chet subdifferential} of the perturbation map in~\cite{CY13}. The \textit{normal subdifferential} of a set-valued map with values in a partially ordered Banach space was introduced in~\cite{Bao_Mordukhovich} by using the Mordukhovich coderivative of the associated \textit{epigraphical multifunction}. If the Mordukhovich coderivative in the definition of normal subdifferential is replaced by the Fr\'{e}chet one, then we have the concept of Fr\'{e}chet subdifferential of a set-valued map. In~\cite{CY13}, Chuong and Yao obtained formulae for outer/inner estimating the Fr\'{e}chet  subdifferential of the perturbation map. Namely, in order to get these results, the \textit{order semicontinuity}  of the objective map and two conditions related to the \textit{solution map} are required. In this paper, due to the convexity assumptions, we get the exact formulas for computing the Fr\'{e}chet  subdifferential of the perturbation map and the conditions obtained are fairly simple if compared to those in~\cite{CY13}. Moreover, our conditions are stated directly on the data of the problem. In addition, our approach is based on convex analysis tools which is different from the approach in~\cite{CY13}.

 The contents of the paper are as follows. Section 2 collects some basic notations and concepts. In Section 3, some calculus rules for the coderivative of convex set-valued maps are stated. It is proved by using the tools from  convex analysis. The main results are presented in Section 4. The coderivative of the perturbation map in vector optimization problems with special constraints are established in Section~5. Some examples are simultaneously provided for analyzing and illustrating the obtained results. In the last section, we give several concluding remarks.

\section{Basic Definitions and Preliminaries}
\markboth{\centerline{\it Sensitivity Analysis  in Parametric  Convex
		Vector Optimization}}{\centerline{\it D.T.V.~An and L.T.~Tung}} \setcounter{equation}{0}
Let $X, Y$ and $P$ be Banach spaces with the duals denoted, respectively, by $X^*$, $Y^*$ and $P^*$. The spaces $X$ and $Y$ are partially ordered corresponding by pointed, closed, and convex cones with apex at the origin $Q\subseteq X$  and $K\subseteq Y$. For $\Omega \subseteq  X$, ${\rm int }\,\Omega,\,{\rm cl} \,\Omega,\, \partial \,\Omega$, and cone\,$\Omega$ denote corresponding the interior, closure, boundary and the cone $\{\lambda a\;|\; \lambda \ge 0,\; a\in \Omega\}$ of $\Omega$.

The polar cone (resp., negative polar cone) $K^*$ (resp., $K^-$) of $K$ is
$$K^*=\{y^*\in Y^*\mid \langle y^*,k\rangle\ge 0,\forall k\in K\},$$
$$(\mbox{resp.,} \ K^-=\{y^*\in Y^*\mid \langle y^*,k\rangle\le 0,\forall k\in K\}).$$
 We write $\mathbb{R}^+(\Omega):=\{ta \in X \mid t\in \mathbb{R}_+, a\in \Omega\}={\rm cone}\,\Omega$, where $\mathbb{R}_+$ signifies the collection of positive numbers.
\medskip

 We recall some concepts of optimality/efficiency, see, e.g., \cite[Chapter 2]{Luc_book_89}, in vector optimization as follows.
\begin{definition} Let $A$ be a nonempty subset of $Y$ and $\bar{a} \in A$. 
\begin{enumerate} \item[(i)] $\bar{a}$ is called a \textit{local (Pareto) minimal/efficient point} of $A$ with respect to (shortly wrt.) $K$, and denoted by  $\bar{a}\in{\rm Min}_KA$, iff there exists $U\in \mathcal{U}(\bar{a})$ such that
$$ (A\cap U-\bar{a})\cap(-K \setminus \{0\})= \emptyset. $$

\item[(ii)] $\bar{a}$ is said to be a \textit{local weak minimal/efficient point} of $A$ wrt. $K$, denoted by $a\in{\rm
    WMin}_KA$, iff there exists $U\in \mathcal{U}(\bar{a})$ such that $$(A\cap U-\bar{a})\cap(-{\rm int}\,K)= \emptyset.$$
\item[(iii)] $\bar{a}$ is termed a \textit{local (Henig) proper minimal/efficient point} of $A$, denoted by $\bar{a} \in {\rm PrMin}_KA$, iff there exists a convex cone
   $Q\subsetneqq Y$ with $K \setminus \{0\}\subseteq {\rm int}\,Q$ and $U\in \mathcal{U}(\bar{a})$  such that
   $$(A\cap U-\bar{a})\cap(-Q\setminus l(Q))= \emptyset,$$
   where $l(Q)=Q \cap (-Q).$
\end{enumerate}
If $U = Y$, the word ``local" is omitted, i.e., we have the corresponding global notions. Throughout the paper, we will focus on global concepts.
\end{definition} 

We have the relationship between above concepts as follows (see~\cite[Proposition 2.2]{Luc_book_89}).
$$ {\rm PrMin}_KA \subset {\rm Min}_KA \subset {\rm
	WMin}_KA.$$

\medskip
In what follows, we will deal with sensitivity analysis of  parameterized vector optimization problems. Firstly, some notations and definitions are recollected. Given  a vector function $f:P\times X\to Y$ and  a set-valued map $C:P\rightrightarrows X$. Let $F:P\rightrightarrows Y$ be a set-valued map defined by
\begin{align}\label{objective map}
	F(p):=F(p,C(p))=\{y\in Y\mid \exists x\in C(p),y=f(p,x)\}.
\end{align}

We consider the following \textit{parametric vector optimization problem}
$$({\rm PSVO}_p)\quad\quad {\rm Min}_K\{ y\in Y\mid \exists x\in C(p),y= f(p,x)\}={\rm Min}_KF(p),$$
\noindent  where $x$ is a decision variable, $p$ is a perturbation parameter, $F$ is an objective map, and $C$ is a constraint map. The \textit{perturbation/frontier map} $\mathcal{F}$, the \textit{weak perturbation/frontier map} $\mathcal{W}$, and the \textit{proper perturbation/frontier map} $\mathcal{P}$ of a family of parameterized vector optimization problem are defined by
$$\mathcal{F}(p):={\rm Min}_KF(p); \ \mathcal{W}(p):={\rm WMin}_KF(p),\;\; \mathcal{P}(p):={\rm PrMin}_KF(p).$$
The \textit{solution map}/the \textit{weak solution map}/the \textit{proper solution map} $\mathcal{S}/\mathcal{S}^w/\mathcal{S}^{\rho}:P\rightrightarrows X$ are given by
$$\mathcal{S}(p):=\{x\in X\mid x\in C(p), f(p,x)\in \mathcal{F}(p)\},$$
$$\mathcal{S}^w(p):=\{x\in X\mid x\in C(p), f(p,x)\in \mathcal{W}(p)\},$$
$$\mathcal{S}^{\rho}(p):=\{x\in X\mid x\in C(p), f(p,x)\in \mathcal{P}(p)\}.$$
\begin{definition} (see, e.g., \cite[Definition 3.1]{Kuk_Tanino_Tanaka_96}, \cite{Tung2016})
For $\bar{p}\in P$ and  a closed convex cone $\widetilde{K}$ contained in ${\rm int} \, K \cup \{0\}.$
\begin{itemize}
	\item [{\rm (i)}] $F$ is said that $K$-\textit{dominated} by $\mathcal{F}$ near $\bar{p}$ iff $F(p)\subseteq \mathcal{F}(p)+K$, for all $p$ in some  neighborhood $\mathcal {U}(\bar{p})$ of $\bar p$.

	\item [{\rm (ii)}] $F$ is said that $\widetilde{K}$-\textit{dominated} by $\mathcal{W}$ near $\bar{p}$ iff $F(p)\subseteq \mathcal{W}(p)+K$, for all $p$ in $\mathcal {U}(\bar{p})$.

	\item [{\rm (iii)}] $F$ is $K$-\textit{dominated} by $\mathcal{P}$ near $\bar{p}$ iff $F(p)\subseteq \mathcal{P}(p)+K$, for all $p$ in $\mathcal {U}(\bar{p})$.
\end{itemize}
\end{definition}

\medskip

The \textit{contingent} (or Bouligand) cone of $\Omega\subset X$ at $\bar{x}\in{\rm cl}\,\Omega$ is
$$T(\bar{x}, \Omega):=\{x\in X\mid \exists\tau_k\downarrow 0,\exists x_k\to x,\;\forall k\in\mathbb{N},\bar{x}+\tau_kx_k\in \Omega\}.$$
For $\bar x\in \Omega$, the \textit{normal cone }to a convex set $\Omega$ at $\bar x$ is
$$N(\bar x, \Omega):=\{ x^*\in X^* \mid \langle x^*, x-\bar x \rangle \le 0 , \ \forall x\in \Omega\}.$$
 In this case, $T(\bar{x},\Omega)={\rm cl\,cone}(\Omega-\bar{x})$, i.e., $T(\bar{x},\Omega)$ is a closed convex cone, and hence,
$N(\bar{x},\Omega)=T(\bar{x},\Omega)^-$ (see, e.g.,~\cite[p. 48]{Bonnans_Shapiro_2000}).

\medskip

Let  $H : X \rightrightarrows Y$ be a set-valued map. The \textit{domain}, \textit{graph}, and \textit{epigraph} of $H$ are given respectively, by
$${\rm dom}\,H=\{x\in X \mid H(x)\not= \emptyset\}, $$
$$ {\rm gph}\,H=\{(x,y)\in X \times Y \mid y \in H(x)\},$$
and
$${\rm epi}\,H:=\{(x,y)\in X\times Y\mid y\in H(x)+K\}.$$
We always assume that $H$ is \textit{proper}, i.e., $\dom H \not= \emptyset.$

\begin{definition} (see, e.g., \cite[Definition 2.1.1]{Aubin_Frankowska_1990}\cite[Chapter 1]{Luc_book_89}\cite[Definition 2.2.4]{Sawaragi_Nakayama_Tanino_85}) Let $H : X \rightrightarrows Y$ be a set-valued map.

\begin{enumerate}
\item[(i)] $H$ is \textit{convex} if ${\rm gph}\,H$ is a convex set in $X\times Y$, i.e., for all $x, x' \in X$ and $\lambda \in [0, 1]$,
    $$\lambda H(x)+(1-\lambda)H(x') \subset H(\lambda x + (1 - \lambda)x').$$
\item[(ii)] $H$ is $K$-\textit{convex }if $\epi{H}$ is convex, i.e., for all $x, x' \in X$ and $\lambda \in [0, 1]$,
		$$ \lambda H(x) + (1 - \lambda)H(x') \subset H(\lambda x + (1 - \lambda)x') + K. $$
\item[(iii)]	$H$ is \textit{closed }(resp., $K$-\textit{closed}) if ${\rm gph}\,H$ (resp., $\epi H$) is a closed set in $X\times Y$.
\end{enumerate}
\end{definition}

	The following proposition can be found in \cite[Proposition 2.1]{Tanino_1988b}, where $X, Y$ and $P$ are finite-dimensional spaces. The same proof remains valid for an infinite-dimensional space setting.
\begin{proposition} \label{Objective_map_convex}
	Suppose that $f:P \times X \rightarrow Y$ is $K$-convex and $C:P \rightrightarrows X$ is a convex set-valued map. Then, $F: P \rightrightarrows Y$ given in~\eqref{objective map} is $K$-convex, i.e., for any $p_1,p_2 \in P$, $\lambda \in [0,1]$,
	$$\lambda F(p_1)+(1-\lambda)F(p_2) \subset F( \lambda p_1 +(1-\lambda)p_2)+K.$$
	In other words, $\gph (F+K)=\epi F$ is convex.
\end{proposition}	

It is worth mentioning here that under the assumptions of Proposition~\ref{Objective_map_convex}, the perturbation map $\mathcal{F}$  can be  nonconvex in general. Let us see some examples.

\begin{example}({\it $f$ is $K$-convex, $C$ is convex and $\mathcal{F}$ is not convex})
	
	Let $X=P=\mathbb{R},$ $Y=\mathbb{R}^2$, $K=\mathbb{R}_+^2$, $f(p,x)=(f_1(p,x),f_2(p,x))=(|x|,|x|),$ and $$C(p)=\{x\in X\mid x\ge |p|\},\ \forall p\in P.$$
	Then, we have
	$$F(p)=\{ y\in Y\mid \exists x\in C(p),y=f(p,x)\}=\{y\in Y\mid y_1=y_2=|x|, x\ge |p|\},$$
	$$\mathcal{F}(p)={\rm Min}_KF(p)=\{y\in Y\mid y_1=y_2=|x|, x=|p|\}=\{(|p|,|p|)\}.$$
	Hence, on one hand for  $p=1,$ $p'=-1\in P$ and $\lambda=\frac{1}{2}\in [0,1]$, one has
	\begin{eqnarray*}
		\lambda \mathcal{F} (p)+(1-\lambda)\mathcal{F} (p')&=&\lambda(|p|,|p|)+
		(1-\lambda)(|p'|,|p'|)\\
		&=& \frac{1}{2}(1,1)+\frac{1}{2}(1,1)=(1,1).
	\end{eqnarray*}
	On the other hand, $ \mathcal{F}(\lambda p+(1-\lambda)p')=(|\lambda p+(1-\lambda)p'|,|\lambda p+(1-\lambda)p'|)=(0,0).$ So $	\lambda \mathcal{F} (p)+(1-\lambda)\mathcal{F} (p') \not\subset  \mathcal{F}(\lambda p+(1-\lambda)p')$.
	Therefore, $\mathcal{F}$ is not convex.
\end{example}
\begin{example}({\it $f$ is $K$-convex, $C$ is convex and $\mathcal{F}$ is convex})
	
	Let $P=\mathbb{R}^3,$ $X=\mathbb{R},$ $Y=\mathbb{R}^2$, $K=\mathbb{R}_+^2$, $$f(p,x)=(f_1(p,x),f_2(p,x))=(x,2x),$$ $$C(p)=\{x\in X\mid x\ge p_1+2p_2+p_3\},\forall p\in P.$$
	Then, we have
	\begin{align*}F(p)&=\{ y\in Y\mid \exists x\in C(p),y=f(p,x)\}\\
		&=\{y\in Y\mid y_1\ge p_1+2p_2+p_3,y_2\ge 2p_1+4p_2+2p_3\},\end{align*}
	$$\mathcal{F}(p)={\rm Min}_KF(p)=\{y\in Y\mid y=(p_1+2p_2+p_3,2p_1+4p_2+2p_3)\},$$
	$$=\{(p_1+2p_2+p_3,2p_1+4p_2+2p_3)\}=\left[\begin{array}{lll}
		1&2&1\\ 2&4&2
	\end{array}
	\right].\left[\begin{array}{l}
		p_1\\ p_2\\p_3
	\end{array}
	\right].$$
	Since $\mathcal{F}$ is a linear map from $\mathbb{R}^3$ to $\mathbb{R}^2$,
	it follows	that $\mathcal{F}$ is convex.
\end{example}

The \textit{profile map} of $H:X \rightrightarrows Y$ is $H + K:X \rightrightarrows Y$ and defined by $(H + K)(x) := H(x) + K$. We have $\gph (H+K)=\epi H$. The following properties follow directly from the definition.
\begin{proposition}\label{convex_K_convex}
	\begin{enumerate}
		\item[\rm (i)] 	The set-valued map $H$ is $K$-convex if and only if $H+K$ is convex.
		\item[\rm (ii)] 	If the set-valued map $H$ is convex, then $H$ is $K$-convex.
	\end{enumerate}
\end{proposition}

Noting that there do exist maps which are $K$-convex but  they are not convex.

\medskip

Next, we recall some concepts from variational analysis.

\begin{definition} {\rm(see \cite[Vol.~I, p. 4 ]{Mordukhovich_2006})}
	\rm
	Let $\Omega$ be a nonempty subset of $X.$ For any $x \in \Omega$, the \textit{Fr\'echet normal cone} of $\Omega$ at $x$ is defined by
	\begin{align*}
		\widehat N(x,\Omega):=\Big\{ x^*\in X^*\mid \limsup\limits_{u \xrightarrow{\Omega}x} \dfrac{\langle x^*, u-x \rangle}{||u-x||} \leq 0 \Big\}.
	\end{align*}
	 If $x \not\in \Omega$, we put$\widehat N(x,\Omega)=\emptyset$.
\end{definition}

Equipping the product space $X\times Y$ with the norm $\|(x,y)\|:=\|x\|+\|y\|$, by the above notion of Fr\'echet normal cone one can define the concept of Fr\'{e}chet coderivative of set-valued maps as follows.

\begin{definition}{\rm (see \cite[Vol. I, p. 40, 41]{Mordukhovich_2006})
	 The \textit{Fr\'{e}chet coderivative} of $H$ at $(\bar x, \bar y) \in {\rm{gph}}\, H$ is the multifunction $\widehat D^* H(\bar x, \bar y): Y^* \rightrightarrows X^*$ given by
		\begin{align*}
			\widehat D^* H(\bar x, \bar y)(y^*):=\left\{x^* \in X^* \mid (x^*, -y^*) \in \widehat N ( (\bar x, \bar y),{\rm{gph}}\, H)\right\}, \ \forall y^* \in Y^*.
		\end{align*}
We omit $\bar y$ in the coderivative notation if $H (x)$ is a singleton.}
\end{definition}

We end this section with the definition of the coderivative for a convex set-valued map. The coderivative concept has not been investigated in the standard framework of convex analysis while it has been well recognized in extended settings of variational analysis; see, e.g., \cite{Mordukhovich_2006},\cite{Rockafellar_Wets_2004}.
\begin{definition}
 (see \cite{An_Yen_2015},\cite[Definition 3.19]{Mordukhovich_Nam_2022}) Let $H : X \rightrightarrows Y$ be a proper convex set-valued map and $(\bar{x},\bar{y})\in {\rm gph}\,H$.
The
coderivative of the $H$ at $(\bar{x},\bar{y})$ is defined by
$$D^*H(\bar{x},\bar{y})(y^*):=\{x^*\in X^*\mid (x^*,-y^*)\in N((\bar{x},\bar{y}),{\rm gph}\,H)\}, \forall y^*\in Y^*,$$
where $N((\bar{x},\bar{y}),{\rm gph}\,H)$ stands for the normal cone of the convex set ${\rm gph}\,H$ at $(\bar x, \bar y)$,
 i.e.
$$N((\bar{x},\bar{y}),{\rm gph}\,H)\!\!=\!\!\{(x^*,y^*) \in X^*\! \times Y^*\! \mid\! \langle x^*, x-\bar x \rangle + \langle y^*, y-\bar y \rangle \le 0, \forall (x,y)\!\in\! {\rm gph}\,H \}.$$
\end{definition}

Since $\gph (H+K)=\epi H$ from the above definition the Fr\'{e}chet coderivative of $H+K$ at $(\bar{x},\bar{y})\in\gph (H+K)$ is the set
$$\widehat{D}^*(H+K)(\bar{x},\bar{y})(y^*)=\{x^*\in X^*\mid (x^*,-y^*)\in \widehat{N}((\bar{x},\bar{y})),{\rm epi}\,H)\}, \forall y^*\in Y^*.$$
In the case, where $H$ is a convex set-valued map, the Fr\'{e}chet coderivative of $H+K$ coincides with the  coderivative   in the sense of convex analysis, i.e., $\widehat{D}^*(H+K)(\bar{x},\bar{y})(y^*) \equiv D^*(H+K)(\bar{x},\bar{y})(y^*),$ for all $y^*\in Y^*.$
\begin{remark}\label{Remark1}
The Fr\'{e}chet coderivative of $H+K$ is known as the \textit{Fr\'{e}chet subdifferential }of $H$ in~\cite{CY13}. The \textit{normal subdifferential} of a set-valued mapping with values in a partially ordered Banach space was introduced in~\cite{Bao_Mordukhovich} by using the Mordukhovich coderivative of the associated \textit{epigraphical multifunction}. If the Mordukhovich coderivative in the definition of normal subdifferential is replaced by the Fr\'{e}chet one, then we have the concept of Fr\'{e}chet subdifferential of a set-valued map.
	\end{remark}

\section{Coderivative Calculus for Convex Set-Valued Maps}
\markboth{\centerline{\it Sensitivity Analysis  in Parametric  Convex
		Vector Optimization}}{\centerline{\it D.T.V.~An and L.T.~Tung}} \setcounter{equation}{0}
The main purpose of this section is to establish some calculus rules for the coderivative of convex set-valued maps, which are needed in the rest of this paper.

In the following propositions  we recall some results on the sum rule and chain rule for coderivative from \cite{MNR17} (see also in~\cite[Chapter 3]{Mordukhovich_Nam_2022}).
\begin{proposition}\label{sum_rule}{\rm (see \cite[Theorem 8.1]{MNR17})} Let $H:X\rightrightarrows Y, $ $L:X\rightrightarrows Y$ be convex set-valued maps between Banach spaces. Suppose that the maps $H$ and $L$ are closed, and the set $\mathbb{R}^+({\rm dom}\,H-{\rm dom}\,L)$ is a closed subspace of $X$.
	Then, for any $(\bar{x},\bar{y})\in {\rm gph}\,(H+L)$ and $y^*\in Y^*$, we have
	$$D^*(H+L)(\bar{x},\bar{y})(y^*)=D^*H(\bar{x},\bar{y}_1)(y^*)
	+D^*L(\bar{x},\bar{y}_2)(y^*),\ (\bar{y}_1,\bar{y}_2)\in S(\bar{x},\bar{y}),$$
	where $S(\bar{x},\bar{y}):=\{(\bar{y}_1,\bar{y}_2)\in Y\times Y\mid \bar y=\bar y_1 +\bar y_2, \ \bar{y}_1\in H(\bar{x}),\bar{y}_2\in L(\bar{x})\}$.
\end{proposition}

Given $H : X \rightrightarrows Y$ and $L: Y\rightrightarrows Z$, the composition $L\circ H:X\rightrightarrows Z$ of $H$ and $L$ is defined by
$$(L\circ H)(x)= \bigcup\limits_{y\in H(x)} L(y):=\{z\in Z \mid \exists y \in H(x) \ \mbox{with}\ z\in L(y) \}, \ x\in X.$$
From the definition, it is easy to see that if $H$ and $L$ are convex, then $L\circ H$ is convex as well.
\begin{proposition}{\rm(see \cite[Theorem 8.2]{MNR17})} \label{chain_rule} Let $H:X\rightrightarrows Y, L:Y\rightrightarrows Z$ be convex set-valued maps between Banach spaces.  Assume that the maps $H$ and $L$ are closed, and the set $\mathbb{R}^+({\rm reg}\,H-{\rm dom}\,L)$ is a closed subspace of $Y$.
	Then, for any $(\bar{x},\bar{z})\in {\rm gph}\,(L\circ H)$ and $z^*\in Z^*$, we have
	$$D^*(L\circ H)(\bar{x},\bar{z})(z^*) =\bigcup\limits_{y^*\in D^*L(\bar{y},\bar{z})(z^*)}D^*H(\bar{x},\bar{y})(y^*),\ \bar{y}\in T(\bar{x},\bar{z}),$$
	where $T(\bar{x},\bar{z}):=H(\bar{x})\cap L^{-1}(\bar{z})$ and ${\rm rge}\,H:=\bigcup\limits_{x\in X}H(x)$.
\end{proposition}

In the sequel, we will use the following result on the coderivative of the map $H$ given by $H(x)=(H_1(x), H_2(x))$, with $H_1:X\rightrightarrows Y_1, H_2:X\rightrightarrows Y_2$ being the maps between Banach spaces. We provide a direct proof of it for the reader’s convenience.

\begin{proposition}\label{sum_rule_H} Let $X, Y_1,Y_2$ be Banach spaces. Suppose that $H_1:X\rightrightarrows~Y_1,$  $H_2:X\rightrightarrows Y_2$ are convex set-valued maps  and $H: X\rightrightarrows Y=Y_1\times Y_2$ is defined by $H(x)=(H_1(x),H_2(x))$. Assume that the maps ${H}_1$ and ${H}_2$ are closed, and the set $\mathbb{R}^+({\rm dom}\,{H}_1-{\rm dom}\,{H}_2)$ is a closed subspace of $X$.
	Then, for any $(\bar{x},\bar{y})\in {\rm gph}\,H$ and $y^*\in Y^*$, we have
	\begin{align}\label{formula1}
		D^*H(\bar{x},\bar{y})(y^*)\!=\!D^*{H_1}(\bar{x},\bar{y}_1)(y^*_1)
		\!+\!D^*{H_2}(\bar{x},\bar{y}_2)(y^*_2),
	\end{align}
	where $y^*=(y^*_1,y^*_2)$  and $\bar y =(\bar{y}_1,\bar{y}_2),  \bar{y}_1\in H_1(\bar{x}), y_2\in H_2(\bar{x}).$
\end{proposition}
\begin{proof} We first show that $H$ is convex. Indeed, for any $(x,y), (x',y')\in {\rm gph}\,H$ and $\lambda\in [0,1]$, one has
	$$y\in H(x)=(H_1(x),H_2(x)), \  y'\in H(x')= (H_1(x'),H_2(x')),$$
	and hence,
\begin{align*}
	\lambda y+(1-\lambda)y'&\in \lambda (H_1(x),H_2(x)) +(1-\lambda)(H_1(x'),H_2(x'))\\
&	=\big(\lambda H_1(x)+(1-\lambda)H_1(x'),\lambda H_2(x)+(1-\lambda)H_2(x')\big)\\
&	\subset(H_1(\lambda x+(1-\lambda)x'),H_2(\lambda x+(1-\lambda)x'))\\
	&=H(\lambda x+(1-\lambda)x'),
\end{align*}
because of the convexity of $H_1$ and $H_2$. Therefore,
	$$\lambda (x,y)+(1-\lambda)(x',y')=(\lambda x+(1-\lambda)x',\lambda y+(1-\lambda)y')\in {\rm gph}\,H,$$
	i.e., $H$ is convex.

	We now apply Proposition~\ref{sum_rule} to the sum $H=\widetilde{H}_1+\widetilde{H}_2$, where
	$\widetilde{H}_1:X\rightrightarrows Y$ and $\widetilde{H}_2:X\rightrightarrows Y$ are defined by
	$$\widetilde{H}_1(x):=(H_1(x),0) \ \mbox{ and}\ \widetilde{H}_2(x):=(0,H_2(x)).$$ It is clear that $\widetilde{H_1}$ and $\widetilde{H_2}$ are convex. On one hand,
	$${\rm dom}\, H_i=\{x\in X \mid H_i(x)\not=\emptyset\}={\rm dom}\,\widetilde{H_i}, \ i=1,2.$$
	On the other hand, \begin{align*}
		{\rm gph}\, \widetilde{H_1}&=\{(x,y_1,y_2)\in X \times Y \mid (y_1,y_2)\in \widetilde{H_1}(x)=(H_1(x),0)\}\\
		&=\{(x,y_1)\in X \times Y_1\times Y_2 \mid y_1\in H_1(x), y_2=0\}={\rm gph}\, {H_1}\times\{0\},
	\end{align*}
	and
	\begin{align*}
		{\rm gph}\, \widetilde{H_2}&=\{(x,y_1,y_2)\in X \times Y \mid (y_1,y_2)\in \widetilde{H_2}(x)=(0,H_2(x))\}\\
		&=\{(x,y_1,y_2)\in X \times Y_1\times Y_2 \mid y_1=0, y_2\in H_2(x)\}=\{0\}\times \gph H_2.
	\end{align*}
	Hence,  $\widetilde{H}_1$ and $\widetilde{H}_2$ are closed due to the closedness of $H_1$ and $H_2$.
	So all the assumptions of Proposition~\ref{sum_rule} are satisfied. Then
	\begin{equation}
		\begin{split}
			D^*H(\bar{x},\bar{y})(y^*)&=D^*H(\bar{x},(\bar{y}_1,0)+(0,\bar{y}_2))(y^*)\\
			\label{formula2}
			&	= D^*\widetilde{H_1}(\bar{x},(\bar{y}_1,0))(y^*)
			+D^*\widetilde{H_2}(\bar{x},(0,\bar{y}_2))(y^*),
		\end{split}
	\end{equation}
	for all $((\bar{y}_1,0),(0,\bar{y}_2))\in \widetilde{S}(\bar{x},\bar{y}).$ 
	In addition, since $Y^*=(Y_1\times Y_2)^*\cong Y_1^*\times Y_2^*$, see e.g., Proposition 1.29 in~\cite{C13}, one has
	$y^*=(y_1^*,y_2^*)\in Y_1^*\times Y_2^*$. So
	\begin{align*}
			x^*\in D^*\widetilde{H_1}(\bar{x},(\bar{y}_1,0))(y^*)&
		\Leftrightarrow(x^*,-y^*)\in N((\bar{x},(\bar{y}_1,0)),{\rm gph}\,\widetilde{H_1}),
	\end{align*}
or, equivalently,
	\begin{align*}\label{formula2-1}
	 &\langle x^*, x-\bar x \rangle -\langle y^*, (y_1,0)-(\bar y_1,0) \rangle\le 0 , \forall (x,y_1,0)\in {\rm gph}\,\widetilde{H_1}\\
		&\Leftrightarrow \langle x^*, x-\bar x \rangle -\langle y^*_1, y_1-\bar y_1 \rangle-\langle y^*_2, 0\rangle\le 0 , \forall (x,y_1)\in {\rm gph}\,H_1\\
		&\Leftrightarrow (x^*,-y^*_1)\in N((\bar{x},\bar{y}_1),{\rm gph}\,H_1)\\
		&\Leftrightarrow x^*\in D^*H_1(\bar{x},\bar{y}_1)(y^*_1),
	\end{align*}
	and hence, $D^*\widetilde{H_1}(\bar{x},(\bar{y}_1,0))(y^*)  =D^*{H_1}(\bar{x},\bar{y}_1)(y^*_1)$. In the same manner we see that $D^*\widetilde{H_2}(\bar{x},(0,\bar{y}_2))(y^*)  =D^*{H_2}(\bar{x},\bar{y}_2)(y^*_2)$. Therefore, we obtain \eqref{formula1} from \eqref{formula2}. The proof is complete.
	\qed
	
\end{proof}

\medskip

	Let $h:X\to Y$ be a single-valued map, we associate $h$ with a scalarization function wrt. some $y^*\in Y^*$, defined by $\langle y^*, h \rangle (x)=\langle y^*, h(x)\rangle$ for all $x\in X.$ The following property of the scalarization function will be use later on.
	\begin{proposition} \label{scalar_function_convex}{\rm (see \cite[Proposition 6.2]{Luc_book_89})} If $K$ is closed, then $h$ is $K$-convex if and only if the composition $\langle y^*, h \rangle$ is a scalar convex function for every $y^*\in K^*.$
	\end{proposition}

We close this section with a result on the relationship between the coderivative of a convex single-valued map and the subdifferential (in the sense of convex analysis) of its scalarization. 
\begin{proposition}\label{scalar_function_coderivative}
	\begin{enumerate}
		\item[\rm (i)]
		Let $h:X\to Y$ be convex. Then,
		$$D^*h(\bar{x})(y^*)=\partial \langle y^*,h\rangle(\bar{x}),\ \forall y^*\in K^*.$$
		\item[\rm (ii)] Let $h:X\to Y$ be $K$-convex. Then,
		$$D^*(h+K)(\bar{x})(y^*)=\partial \langle y^*,h\rangle(\bar{x}),\ \forall y^*\in K^*.$$
			\end{enumerate}
	In particular, if $Y=\mathbb{R}^m$ and $h:X\to \mathbb{R}^m$ is continuously differentiable with $\nabla h(\bar{x})=(\nabla h_1(\bar{x}),...,\nabla h_m(\bar{x}))$, then 
		$$D^*(h+K)(\bar{x})(y^*)=\partial \langle y^*,h\rangle(\bar{x})=\nabla h(\bar{x})^*(y^*),\ \forall y^*\in K^*,$$
		where $\nabla h(\bar{x})^*$ is the adjoint operator of $\nabla h(\bar x)$.
\end{proposition}
\begin{proof} (i) As $h$ is convex, we deduce from Proposition \ref{convex_K_convex} that $h$ is $K$-convex.
	Since $K$ is a closed, convex cone, by Proposition~\ref{scalar_function_convex}, we have $\langle y^*, h \rangle$ is a scalar convex function for every $y^*\in K^*.$ Pick any $x^*\in \partial \langle y^*,h\rangle(\bar{x})$, by the definition, this is equivalent to
	\begin{align*}
		&\langle x^*, x-\bar x \rangle \le \langle y^*, h \rangle (x) -\langle y^*, h \rangle (\bar x) , \ \forall x\in X\\
		&\Leftrightarrow \langle x^*, x-\bar x\rangle - \langle y^*, h(x)-h (\bar x)\rangle \le 0,  \ \forall x\in X\\
		&\Leftrightarrow(x^*, -y^*)\in N((\bar x, h(\bar x)), {\rm gph}\, h)\\
		&\Leftrightarrow x^* \in D^*h(\bar{x})(y^*).
	\end{align*}
	\noindent (ii) $*$ Pick any $x^*\in \partial \langle y^*,h\rangle(\bar{x})$. Keeping in mind that $\langle y^*,h\rangle$ is convex for every $y^*\in K^*$ by Proposition~\ref{scalar_function_convex}. According to the definition of the  subdifferential, one~has
	\begin{align*}
	\langle x^*, x-\bar x \rangle \le \langle y^*, h \rangle (x) -\langle y^*, h \rangle (\bar x) , \ \forall x\in X.	
	\end{align*}
This implies that
\begin{align*}
&\langle x^*, x-\bar x \rangle \le \langle y^*, h(x) \rangle -\langle y^*, h(\bar x) \rangle \le \langle y^*, y \rangle -\langle y^*, h(\bar x) \rangle, \ \forall y\in h(x)+K,\forall x\in X\\
&\Leftrightarrow \langle (x^*,y^*), (x,y)-(\bar{x},h (\bar x)\rangle \le 0,   \forall y\in h(x)+K,\ \forall x\in X\\
&\Leftrightarrow(x^*, -y^*)\in N((\bar x, h(\bar x)), {\rm epi}\, h)\\
&\Leftrightarrow x^* \in D^*(h+K)(\bar{x})(y^*).
\end{align*}
	Conversely, let $x^* \in D^*(h+K)(\bar{x})(y^*)$. Then,
	\begin{align*}
			&(x^*, -y^*)\in N((\bar x, h(\bar x)), {\rm epi}\, h)\\
		&\Leftrightarrow \langle (x^*,y^*), (x,y)-(\bar{x},h (\bar x)\rangle \le 0, \  \forall y\in h(x)+K,\ \forall x\in X.
	\end{align*}
Noting that $h(x)+\{0\} \in h(x)+K.$ The last inequality implies that
	\begin{align*}	
		&\langle (x^*,y^*), (x,h(x))-(\bar{x},h (\bar x)\rangle \le 0,   \forall y\in h(x)+\{0\}\in h(x)+K,\forall x\in X\\
		&\Leftrightarrow\langle x^*, x-\bar x \rangle \le \langle y^*, h \rangle (x) -\langle y^*, h \rangle (\bar x) , \ \forall x\in X\\
		&\Leftrightarrow x^*\in \partial \langle y^*,h\rangle(\bar{x}).
	\end{align*}
	\noindent $*$ Let $y^*=(y^*_1,...,y^*_m)$ then $\langle y^*,h\rangle=\sum_{i=1}^m\limits y^*_ih_i:X\to \mathbb{R}$. Hence, $\langle y^*,h\rangle$ is also continuously differentiable  and
	$\partial \langle y^*,h\rangle(\bar{x})=\sum_{i=1}^m\limits y^*_i\nabla h_i(\bar{x})=\nabla h(\bar{x})^*(y^*).$
	\qed
\end{proof}

\section{The Main Results}
\markboth{\centerline{\it Sensitivity Analysis  in Parametric  Convex
		Vector Optimization}}{\centerline{\it D.T.V.~An and L.T.~Tung}} \setcounter{equation}{0}
In this section we will give formulas for computing the Fr\'{e}chet coderivative of the profile of the perturbation map~$\mathcal{F}$ (resp., the weak perturbation map $\mathcal{W}$, the proper perturbation map $\mathcal{P}$) at the corresponding point.

\medskip

The first theorem provides us with formulas for computing the coderivative of the profile map of $F$. This is the main key to obtaining the next results.
\begin{theorem}\label{Thm_objective_map}
Let $f:P\times X\to Y$ be  a $K$-convex and $K$-closed map, $C:P\rightrightarrows X$ be closed and convex, $(\bar{p},\bar{x})\in {\rm gph}\,\mathcal{S}$ and $\bar{y}=f(\bar{p},\bar{x})$. Let $H:P\rightrightarrows P\times X$ be defined by $H(p)=(H_1(p),H_2(p))$, where $H_1(p)=\{p\}$ and $H_2(p)=C(p)$. Suppose that the following conditions:
\begin{enumerate}
\item[\rm (i)] $\mathbb{R}^+({\rm reg}\,H-{\rm dom}\,f)$ is a closed subspace of $P\times X$,
\item[\rm (ii)] $\mathbb{R}^+(P-{\rm dom}\,C)$ is a closed subspace of $P$
\end{enumerate}
 are satisfied.
Then, for any $y^*\in Y^*$, we have
\begin{equation}\label{prop2-4_new}
	D^*(F+K)(\bar{p},\bar{y})(y^*)\!=\! \bigcup\limits_{(p^*,x^*)\in D^*(f+K)((\bar{p},\bar{x}))(y^*)}
	\bigg\{p^*
	+D^*C(\bar{p},\bar{x})(x^*)\bigg\}.
\end{equation}
In particular, for $y^*\in K^*$, one has
\begin{equation}\label{prop2-1-1}
	D^*(F+K)(\bar{p},\bar{y})(y^*)= \bigcup\limits_{(p^*,x^*)\in \partial \langle y^*,f\rangle (\bar{p},\bar{x})}\bigg\{p^*
+D^*C(\bar{p},\bar{x})(x^*)\bigg\}.
\end{equation}
Moreover, if $f$ is continuously differentiable at $(\bar{p},\bar{x})$, then
$$	D^*(F+K)(\bar{p},\bar{y})(y^*)= \nabla_p f(\bar{p},\bar{x})^*(y^*)+D^*C(\bar{p},\bar{x})\big(\nabla_x f(\bar{p},\bar{x})^*(y^*)\big), \forall y^*\in K^*.$$
\end{theorem}
\begin{proof}
We first observed that $F+K=(f+K)\circ H$. We now check assumptions of Proposition~\ref{chain_rule} to get the formula for coderivative of $F+K$. The convexity and closedness of $f+K$ follow from Proposition~\ref{convex_K_convex} and our assumptions. While by the same manner as in the proof of Proposition~\ref{sum_rule_H}, it follows  that $H$ is a closed convex map. Hence, if condition~(i) is satisfied, we deduce from Proposition~\ref{chain_rule} that, for any $y^*\in Y^*$,
\begin{equation}\label{prop2-3}
D^*(F+K)(\bar{p},\bar{y})(y^*)\!=\! \bigcup\limits_{(p^*,x^*)\in D^*(f\!+\!K)((\bar{p},\bar{x}))(y^*)}
D^*H(\bar{p},(\bar{p},\bar{x}))(p^*,x^*).
\end{equation}
Since $H_1:P\rightrightarrows P$, $H_1(p)=p$ for all $p\in P$, it follows that $H_1\equiv id_P$, where $id_P$ is the identity map, and ${\rm dom}\,H_1=P$. Moreover, we deduce from $H_2:P\rightrightarrows X$ and $H_2(p)=C(p)$ for all $p\in P$ that $H_2\equiv C$.
Hence, $$\mathbb{R}^+({\rm dom}\,{H}_1-{\rm dom}\,{H}_2)\equiv\mathbb{R}^+(P-{\rm dom}\,C).$$ Keeping in mind that $H_1$, $H_2$ are closed convex, combining  the latter with the validity of condition~(ii), we see  that all the assumptions of  Proposition~\ref{sum_rule_H} are satisfied. Hence,
\begin{equation*}
	D^*H(\bar{p},(\bar{p},\bar{x}))(p^*,x^*)=D^*H_1(\bar{p},\bar{p})(p^*)+
	D^*H_2(\bar{p},\bar{x})(x^*).
\end{equation*}
By~\cite[Lemma 49]{Li_Penot_Xue}, we have $ D^*H_1(\bar{p},\bar{p})(p^*)\!+\!
D^*H_2(\bar{p},\bar{x})(x^*)\!=\!p^*\!+\!D^*C(\bar{p},\bar{x})(x^*).$
Consequently, 
\begin{equation}\label{prop2-2}
	D^*H(\bar{p},(\bar{p},\bar{x}))(p^*,x^*)=p^*+D^*C(\bar{p},\bar{x})(x^*).
\end{equation}
Combining \eqref{prop2-3} and \eqref{prop2-2} we obtain~\eqref{prop2-4_new}.

Since $f$ is $K$-convex, by invoking Proposition \ref{scalar_function_coderivative} (ii), one has
$$ D^*(f+K)((\bar{p},\bar{x}))(y^*)
=\partial^*f((\bar{p},\bar{x}))(y^*)=\partial \langle y^*,f\rangle (\bar{p},\bar{x}), \ \forall y^*\in K^*.$$
Combining the latter with~\eqref{prop2-4_new} we get~\eqref{prop2-1-1}.

In particular, if $f$ is continuously differentiable at $(\bar{p},\bar{x})$ then
$$\partial\langle y^*,f\rangle (\bar{p},\bar{x})=\nabla f(\bar{p},\bar{x})^*(y^*)=\big(\nabla_pf(\bar{p},\bar{x})^*(y^*)
,\nabla_xf(\bar{p},\bar{x})^*(y^*)\big).$$
Therefore the conclusion is obtained.
	\qed
\end{proof}
\medskip

We now establish formulae for computing the  Fr\'{e}chet coderivative of the profile of the perturbation map~$\mathcal{F}$.
\begin{theorem}\label{Thm_frontier_map} In addition to all the assumptions of Theorem~\ref{Thm_objective_map}, suppose
that $F$ is $K$-dominated by $\mathcal{F}$ near $\bar{p}$.
Then, for any $y^*\in Y^*$, we have
\begin{equation}\label{2.3-1_new}
\widehat{D}^*(\mathcal{F}+K)(\bar{p},\bar{y})(y^*)\!=\!\bigcup\limits_{(p^*,x^*)\in D^*(f+K)((\bar{p},\bar{x}))(y^*)}
	\bigg\{p^*
	+D^*C(\bar{p},\bar{x})(x^*)\bigg\}.
\end{equation}
In particular, for all $y^*\in K^*$, one has
\begin{equation}\label{2.3-1}
\widehat{D}^*(\mathcal{F}+K)(\bar{p},\bar{y})(y^*)= \bigcup\limits_{(p^*,x^*)\in \partial \langle y^*,f\rangle (\bar{p},\bar{x})}\bigg\{p^*
+D^*C(\bar{p},\bar{x})(x^*)\bigg\}.
\end{equation}
Moreover, if $f$ is continuously differentiable at $(\bar{p},\bar{x})$, then
$$\widehat{D}^*(\mathcal{F}+K)(\bar{p},\bar{y})(y^*)= \nabla_p f(\bar{p},\bar{x})^*(y^*)+D^*C(\bar{p},\bar{x})\big(\nabla_x f(\bar{p},\bar{x})^*(y^*)\big),\  \forall y^*\in K^*.$$
\end{theorem}
\begin{proof} Since $\mathcal{F}(p)\subset F(p)$ for all $p\in P$ and the domination property holds, there
exists a neighborhood $U$ of $\bar{p}$ such that
$$\mathcal{F}(p)+K=F(p)+K, \ \forall p\in U.$$
Hence
$$\widehat{D}^*(\mathcal{F}+K)(\bar{p},\bar{y})(y^*)
=\widehat{D}^*(F+K)(\bar{p},\bar{y})(y^*),\ \forall y^*\in Y^*.$$
Under the assumptions, by Proposition~\ref{Objective_map_convex}, $F$ is $K$-convex. Moreover, thanks to Proposition~\ref{convex_K_convex}, we have $F+K$ is convex. So,
\begin{equation}\label{2.3-2}
\widehat{D}^*(F+K)(\bar{p},\bar{y})(y^*)
\equiv D^*(F+K)(\bar{p},\bar{y})(y^*),\ \forall y^*\in Y^*,
\end{equation}
and hence, we obtain~\eqref{2.3-1_new} from~\eqref{2.3-2} and~\eqref{prop2-4_new}.
In the case, where $y^*\in K^*$, we get from Proposition~\ref{scalar_function_coderivative} that
$$D^*(f+K)((\bar{p},\bar{x}))(y^*)=\partial \langle y^*,f\rangle (\bar{p},\bar{x}).$$
Therefore, combining~\eqref{prop2-1-1} and~\eqref{2.3-2} yields~\eqref{2.3-1}. The proof is complete.
	\qed
\end{proof}

\begin{remark} In~\cite{CY13},  the authors call $\widehat{D}^*(\mathcal{F}+K)(\bar{p},\bar{y})$ at the point under consideration  the\textit{ Fr\'{e}chet  subdifferential} as mentioned in Remark~\ref{Remark1}. By using some suitable conditions, the authors obtained  outer/inner estimations for the Fr\'{e}chet  subdifferential of the perturbation map~$\mathcal{F}$. Namely, the \textit{order semicontinuity}  of $F$ and two conditions related to the solution map $\mathcal{S}$ are required. In this paper, due to the convexity assumptions, we get the exact formulas for computing the Fr\'{e}chet  subdifferential of the perturbation map and the conditions set are fairly simple if compared to those in~\cite{CY13}. Moreover, our conditions are stated directly on the data of the problem. In addition, our approach is based on convex analysis tools which is  different from the approach in~\cite{CY13}.
\end{remark}

\medskip
By some suitable changes,  the results in this paper still hold for the weak perturbation/solution map or the proper perturbation/solution map.
\begin{theorem} \label{Weak_theorem} In addition to all the assumptions of Theorem~\ref{Thm_objective_map}, suppose
that $F$ is $\widetilde{K}$-dominated by $\mathcal{W}$ near $\bar{p}$.
Then, for any $y^*\in Y^*$, we have
\begin{equation*}
\widehat{D}^*(\mathcal{W}+ \tilde {K})(\bar{p},\bar{y})(y^*)\!=\!\bigcup\limits_{(p^*,x^*)\in D^*(f+\tilde {K})((\bar{p},\bar{x}))(y^*)}
	\bigg\{p^*
	+D^*C(\bar{p},\bar{x})(x^*)\bigg\}.
\end{equation*}
In particular, for any $y^*\in  \tilde {K}^*,$
\begin{equation*}
\widehat{D}^*(\mathcal{W}+ \tilde {K})(\bar{p},\bar{y})(y^*)= \bigcup\limits_{(p^*,x^*)\in \partial \langle y^*,f\rangle (\bar{p},\bar{x})}\bigg\{p^*
+D^*C(\bar{p},\bar{x})(x^*)\bigg\}, \forall y^*\in  \tilde {K}^*.
\end{equation*}
Moreover, $f$ is continuously differentiable at $(\bar{p},\bar{x})$, then
$$\widehat{D}^*(\mathcal{W}+ \tilde {K})(\bar{p},\bar{y})(y^*)= (\nabla_p f(\bar{p},\bar{x})^*)(y^*)+D^*C(\bar{p},\bar{x})(\nabla_x f(\bar{p},\bar{x})^*(y^*)), \forall y^*\in \tilde {K}^*.$$
\end{theorem}
\begin{proof}
For $(\bar{p},\bar{x})\in {\rm gph}\,\mathcal{S}^w$ and $\bar y=f(\bar p,\bar x)$.  If $F$ is $\widetilde{K}$-dominated by $\mathcal{W}$ near $\bar{p}$ and $\mathcal{W}+\tilde K$ is convex, we have $$\widehat{D}^*(\mathcal{W}+\tilde K)(\bar{p},\bar{y})(y^*)
		=\widehat{D}^*(F+\tilde K)(\bar{p},\bar{y})(y^*)$$ for any  $ y^*\in Y^*$.
Thus, by Theorem~\ref{Thm_objective_map}, we get desired results.
	\qed
	\end{proof}
\begin{theorem} \label{Proper_theorem} In addition to all the assumptions of Theorem~\ref{Thm_objective_map}, suppose
that $F$ is $K$-dominated by $\mathcal{P}$ near $\bar{p}$.
Then, for any $y^*\in Y^*$, we have
\begin{equation*}
	\widehat{D}^*(\mathcal{P}+  K)(\bar{p},\bar{y})(y^*)\!=\!\bigcup\limits_{(p^*,x^*)\in D^*(f+K)((\bar{p},\bar{x}))(y^*)}
	\bigg\{p^*
	+D^*C(\bar{p},\bar{x})(x^*)\bigg\}.
\end{equation*}
In particular, for any $y^*\in  K^*,$
\begin{equation*}
	\widehat{D}^*(\mathcal{P}+  K)(\bar{p},\bar{y})(y^*)= \bigcup\limits_{(p^*,x^*)\in \partial \langle y^*,f\rangle (\bar{p},\bar{x})}\bigg\{p^*
+D^*C(\bar{p},\bar{x})(x^*)\bigg\}, \forall y^*\in K^*.
\end{equation*}
In particular, if $f$ is continuously differentiable at $(\bar{p},\bar{x})$, then
$$	\widehat{D}^*(\mathcal{P}\!+\!  K)(\bar{p},\bar{y})(y^*)\!=\! (\nabla_p f(\bar{p},\bar{x})^*)(y^*)\!+\!D^*C(\bar{p},\bar{x})(\nabla_x f(\bar{p},\bar{x})^*(y^*)), \forall y^*\in K^*.$$
\end{theorem}
\begin{proof}
For $(\bar{p},\bar{x})\in {\rm gph}\,\mathcal{S}^\rho$ and $\bar y=f(\bar p,\bar x)$.   If $F$ is ${K}$-dominated by $\mathcal{P}$ near $\bar{p}$ and $\mathcal{P}+K$ is convex, we have $$\widehat{D}^*(\mathcal{P}+K)(\bar{p},\bar{y})(y^*)
	=\widehat{D}^*(F+K)(\bar{p},\bar{y})(y^*)$$ for any $ y^*\in Y^*$.
	Hence, the proof is immediate from Theorem~\ref{Thm_objective_map}.
		\qed
\end{proof}
\medskip

We now give an example, which aims to illustrate Theorem~\ref{Thm_frontier_map}.
\begin{example}
Let $P=\mathbb{R}^3,X=\mathbb{R},Y=\mathbb{R}^2$, $Q=\mathbb{R}_+,K=\mathbb{R}_+^2$, $$f(p,x)=(f_1(p,x),f_2(p,x))=(x,2x),$$ $$C(p)=\{x\in X\mid x\ge p_1+2p_2+p_3\},\forall p\in P.$$
Then, we have
\begin{align*}
	F(p)&=\{ y\in Y\mid \exists x\in C(p),y=f(p,x)\}\\
	&=\{y\in Y\mid y_1\ge p_1+2p_2+p_3,y_2\ge 2p_1+4p_2+2p_3\},
\end{align*}
$$\mathcal{F}(p)={\rm Min}_KF(p)=\{y\in Y\mid y=(p_1+2p_2+p_3,2p_1+4p_2+2p_3)\},$$
$$=\{(p_1+2p_2+p_3,2p_1+4p_2+2p_3)\}=\left[\begin{array}{lll}
1&2&1\\ 2&4&2
\end{array}
\right].\left[\begin{array}{l}
p_1\\ p_2\\p_3
\end{array}
\right],$$
\begin{align*}
	\mathcal{S}(p)&=\{x\in X\mid x\in C(p), f(p,x)\in \mathcal{F}(p)\}\\
	&=\{x\in X\mid x\ge p_1+2p_2+p_3,(x,2x)=(p_1+2p_2+p_3,2p_1+4p_2+2p_3)\}\\
	&=\{p_1+2p_2+p_3\},
\end{align*}
and hence, the domination property holds for $F$ around $\bar{p}$. Moreover, we deduce from the fact that $\mathcal{F}$ is a linear map from $\mathbb{R}^2$ to $\mathbb{R}^2$
that $\mathcal{F}$ is convex.

Taking $\bar{p}=(0,0,{0})$ and $\bar{x}=0\in \mathcal{S}(\bar{p})=\{0\}$, then $\bar{y}=f(\bar{p},\bar{x})=(0,0)$ one~has
\begin{equation}\label{Coderivative}
	\begin{split}
D^*\mathcal{F}(\bar{p},\bar{y})(y^*)&=
\left[\begin{array}{lll}
1&2&1\\ 2&4&2
\end{array}
\right]^T.\left[\begin{array}{l}
y_1^*\\ y_2^*
\end{array}
\right]=\left[\begin{array}{ll}
1&2\\ 2&4\\ 1&2
\end{array}
\right].\left[\begin{array}{l}
y_1^*\\ y_2^*
\end{array}
\right]\\
&=(y_1^*+2y_2^*,2y_1^*+4y_2^*,y_1^*+2y_2^*),
\end{split}
\end{equation} 
and
$\partial f_1(\bar{p},\bar{x})=\{((0,0,0),1)\}, \  \partial f_2(\bar{p},\bar{x})=\{((0,0,0),2)\}.$
Since $C$ is a convex
set-valued mapping, we have
$$N((\bar{p},\bar{x}),{\rm gph}\,C)=\{(1,2,1,-1)x^*\mid x^*\ge 0\}$$
and
\begin{align*}D^*C(\bar{p},\bar{x})(x^*)&=\{p^*\in P^*\mid (p^*,-x^*)\in N((\bar{p},\bar{x},{\rm gph}\,C)\}\\
&=\left\{
\begin{array}{ll}
\{(x^*,2x^*,x^*)\}, \ \mbox{if}\; x^*\ge 0,\\
\emptyset, \ \mbox{otherwise.}
\end{array}
\right.
\end{align*}
We now check the assumption (i) and (ii) of Theorem~\ref{Thm_objective_map}. On one hand, as
${\rm reg}H=\bigcup\limits_{q\in P}H(q)=\{(p,x)\in P\times X\mid p=q,x\ge p_1+2p_2+p_3\},$
and ${\rm dom}\,f=\mathbb{R}^4,$
we see that $\mathbb{R}^+({\rm reg}H-{\rm dom}\,f)=\mathbb{R}^4$ is a closed subspace of $P\times X$.
On the other hand, since ${\rm dom}\,C=P=\mathbb{R}^3$, one has $\mathbb{R}^+(P-{\rm dom}\,C)=P$ is a closed subspace of $P$. Therefore, all assumptions of
Theorem~\ref{Thm_objective_map} are satisfied.

 Let $y^*=(y^*_1,y^*_2)\in K^*=\mathbb{R}^2_+$. Then, $y^*_1\ge 0, y^*_2\ge 0$ and $$\langle y^*,f\rangle(p,x) = y^*_1f_1(p,x)+y^*_2f_2(p,x)=y^*_1x+y^*_22x=(y^*_1+2y^*_2)x,$$
$$\partial \langle y^*,f\rangle(\bar{p},\bar{x})=\{((0,0,0),y^*_1+2y^*_2)\},$$
\begin{equation*}
	\begin{split}
\bigcup\limits_{(p^*,x^*)\in \partial \langle y^*,f\rangle(\bar{p},\bar{x})}\{p^*\!+\!D^*C(\bar{p},\bar{x})(x^*)\}&=
(0,0,0)\!+\!D^*C(\bar{p},\bar{x})(y^*_1+2y^*_2)\\
&=(y^*_1\!+\!2y^*_2,2(y^*_1\!+\!2y^*_2),y^*_1\!+\!2y^*_2).	
\end{split}
\end{equation*} 
Hence, from the latter and~\eqref{Coderivative}, we conclude that
$$D^* \mathcal{F}(\bar{p},\bar{y})(y^*)=\bigcup\limits_{(p^*,x^*)\in \partial \langle y^*,f\rangle(\bar{p},\bar{x})}\{p^*+D^*C(\bar{p},\bar{x})(x^*)\}.$$

\end{example}

\section{Applications}
\markboth{\centerline{\it Sensitivity Analysis  in Parametric  Convex
		Vector Optimization}}{\centerline{\it D.T.V.~An and L.T.~Tung}} \setcounter{equation}{0}
\subsection{Constraints Described by Finitely Many Equalities and Inequalities}
We consider problem (${\rm PSVO}_p)$ with the functional constraints described by finitely many equalities and inequalities given as follows
\begin{equation}\label{constraint_C}
	\begin{split}C(p)=\{x\in X\mid &g_i(p,x)\le 0,\  i=1,..., m, \\
	& g_i(p,x)=0,\ i=m+1,..., m+r \},
\end{split}
\end{equation}
where  $g_i:P\times X\to \mathbb{R}$ are continuously differentiable and convex functions.
It is clear that $C$ is a closed convex set-valued map.

Recall that the set-valued map $C$ satisfies the \textit{Abadie constraint qualification }(ACQ) at $(\bar{p},\bar{x})\in {\rm gph}\,C$ if
\begin{equation}\label{ACQ}
	\begin{split}
T((\bar{p},\bar{x}),{\rm gph}\,C)& \!\supset\! \bigg\{d\!\in\! P\times X\mid \langle \nabla g_i(\bar{p},\bar{x}),d\rangle\le 0, \, \mbox{\rm for}\, i\in I(\bar p, \bar x),\\
&\ \quad \quad \quad
 \langle\nabla g_i(\bar{p},\bar{x}),d\rangle\!=\! 0,\, \mbox{\rm for}\, i\!=\!m\!+\!1,...,m\!+\!r\bigg\},
\end{split}
\end{equation}
where $I(\bar p, \bar x)=\{i=1,...,m \mid g_i(\bar p, \bar x)=0\}.$

 In~\cite{WuLi97}, the author shows that (ACQ) is the weakest condition that ensures the characterization of an optimal solution by Karush-Kuhn-Tucker (KKT) conditions.
 
 \medskip

The following proposition gives us the formula for computing the coderivative of the constraint map $C$ given by~\eqref{constraint_C} under the validity of the (ACQ) condition.
\begin{proposition}\label{prop4.1}
	Let $C$ be given by~\eqref{constraint_C} and $(\bar{p},\bar{x})\in {\rm gph}\, C$. Suppose that $C$ satisfies the Abadie constraint qualification~\eqref{ACQ} at $(\bar{p},\bar{x})$. Then,
\begin{equation*}
	\begin{split}D^*C(\bar{p},\bar{x})(x^*)\!=\!\bigg\{p^*\in P^*&\mid (p^*, -x^*)\!=\!\sum\limits_{i=1}^{m+r}\lambda_i\nabla g_i(\bar{p},\bar{x}) \ \mbox{\rm for some} \ \lambda\in \mathbb{R}^{m+r}\\
		&  \quad \quad \ \mbox{\rm with} \ \lambda_i \ge 0, \lambda_ig_i(\bar p, \bar x)=0 \ \mbox{\rm as}\ i=1,...,m \bigg\}.
	\end{split}
\end{equation*}
\end{proposition}
\begin{proof} Let $d=(p,x)\in {T}((\bar{p},\bar{x}),{\rm gph}\,C)$. Then, by the definition, there {exist} $\tau_k\downarrow 0$ and $(p_k,x_k)\to (p,x)$ such that
$$\bar{x}+t _kx_k\in C(\bar{p}+t _kp_k), \ \forall k\in \mathbb{N}.$$
{The latter means that}
$$g_i(\bar{p}+t _kp_k,\bar{x}+t _kx_k)\le 0,  \ \forall k\in \mathbb{N},\ \forall i=1,..,m,$$
$$g_i(\bar{p}+t _kp_k,\bar{x}+t _kx_k)= 0,  \ \forall k\in \mathbb{N},\ \forall i=m+1,...,m+r.$$
We deduce from the continuously Fr\'{e}chet differentiability of $ g_i$ that
$$g_i(\bar{p},\bar{x})+t_k\nabla g_i(\bar{p},\bar{x})(p_k,x_k)+o(t _k\|(p_k,x_k)\|) \le 0,  \ \forall k\in \mathbb{N},\ \forall i=1,...,m,$$
$$g_i(\bar{p},\bar{x})+t_k\nabla g_i(\bar{p},\bar{x})(p_k,x_k)+o(t _k\|(p_k,x_k)\|) = 0,  \ \forall k\in \mathbb{N},\ \forall i=m+1,...,m+r,$$
where $\lim\limits_{k\to\infty} \dfrac{o(t _k\|(p_k,x_k)\|)}{t_k}=0. $
Since $g_i(\bar{p},\bar{x})=0$ for all $i \in I(\bar{p},\bar{x})$ and $g_i(\bar{p},\bar{x})=0,$ for every $ i=m+1,...,m+r$, one has
$$\nabla g_i(\bar{p},\bar{x})(p_k,x_k)+\frac{1}{t _k}o(t _k\|(p_k,x_k)\|) \le 0, \forall k\in  \mathbb{N}, \forall i\in I(\bar{p},\bar{x}),$$
$$\nabla g_i(\bar{p},\bar{x})(p_k,x_k)+\frac{1}{t _k}o(t _k\|(p_k,x_k)\|) = 0, \forall k\in \mathbb{N}, \forall i=m+1,...,m+r.$$
{Letting $k\to \infty$, the above inequalities and equalities derive that}
$$\nabla g_i(\bar{p},\bar{x})(p,x)\le 0, \forall i\in I(\bar{p},\bar{x}) \ \mbox{and} \ \nabla g_i(\bar{p},\bar{x})(p,x)\le 0,\ \forall i=m+1,..,m+r.$$
Hence,
\begin{align*}T((\bar{p},\bar{x}),{\rm gph}\,C) &\subset\big\{d\in P\times X\mid \langle \nabla g_i(\bar{p},\bar{x}),d\rangle\le 0, \,\mbox{\rm for}\ i\in I(\bar{p},\bar{x}), \\
		& \quad \quad \quad \quad \langle\nabla g_i(\bar{p},\bar{x}),d\rangle= 0, \ \mbox{\rm for}\ i=m+1,...,m+r\big\}.
	\end{align*}
Moreover, since $C$ satisfies the Abadie constraint qualification (ACQ) at $(\bar{p},\bar{x})$, we have
\begin{align*}&T((\bar{p},\bar{x}),{\rm gph}\,C)\\
	& \!=\!\{d\!\in\! P\!\times\! X\!\mid \!\langle \nabla g_i(\bar{p},\bar{x}),d\rangle\!\le\! 0, \forall i\in I(\bar{p},\bar{x}), \langle\nabla g_i(\bar{p},\bar{x}),d\rangle\!=\! 0,i=m+1,...,m+r\}\\
&\!=\!\{d\!\in\! P\!\times\! X\!\mid\! \langle \nabla g_i(\bar{p},\bar{x}),d\rangle\!\le\! 0, \forall i\in I(\bar{p},\bar{x}), \langle \pm \nabla g_i(\bar{p},\bar{x}),d\rangle\!\le\! 0,i=m+1,...,m+r\}.
\end{align*}
We now take any $p^*\in D^*C(\bar{p},\bar{x})(x^*)$. {This is equivalent to
\begin{align*}
	(p^*,-x^*)\in N((\bar{p},\bar{x}),{\rm gph}\,C).	\end{align*}
 As ${\rm gph}\,C$ is convex, we have \begin{align*}
		N((\bar{p},\bar{x}),{\rm gph}\,C)&= T((\bar{p},\bar{x}),{\rm gph}\,C)^-.
	\end{align*}
Consequently,
\begin{align*}
	(p^*,-x^*)&\in T((\bar{p},\bar{x}),{\rm gph}\,C)^-,
\end{align*}
i.e., $\langle (p^*,-x^*), (p,x) \rangle \le 0$ for all $(p,x)\in T((\bar{p},\bar{x}),{\rm gph}\,C).$ In other words, the inequality $\langle (p^*,-x^*), (p,x) \rangle \le 0$ is a consequence of the inequalities system
$\langle \nabla g_i(\bar{p},\bar{x}),d\rangle\le 0, \forall i\in I(\bar{p},\bar{x}), \langle \pm \nabla g_i(\bar{p},\bar{x}),d\rangle\le 0,i=m+1,...,m+r.$
By~\cite[Lemma 1]{Bartl}, there exist $\lambda_i\ge 0$ such that}
$$(p^*,-x^*)=\sum\limits_{i\in I(\bar p, \bar x)}\lambda_i\nabla g_i(\bar{p},\bar{x})+\sum\limits_{i=m+1}^{m+r}\bigg(\lambda_i\nabla g_i(\bar{p},\bar{x})+\lambda_i\big(-\nabla g_i(\bar{p},\bar{x})\big)\bigg).$$
The latter can be rewritten as
$(p^*,-x^*)=\sum\limits_{i=1}^{m+r}\lambda_i\nabla g_i(\bar{p},\bar{x}), $ for some $ \lambda\in \mathbb{R}^{m+r}$  with $\lambda_i \ge 0, \lambda_i g_i(\bar p, \bar x)=0$ as $ i=1,...,m,$
which completes the proof.
	\qed
\end{proof}

Our next theorem provides a formula for computing the Fr\'{e}chet coderivative of the profile of $\mathcal{F}$ via the classical Lagrange multipliers in the case where constraint functions are Fr\'{e}chet differentiable.
\begin{theorem} Let $f:P\times X\to Y$ be $K$-closed and $K$-convex, $C:P\rightrightarrows X$ be a convex closed set-valued map given by~\eqref{constraint_C}, $(\bar{p},\bar{x})\in {\rm gph}\,\mathcal{S}$ and $\bar{y}=f(\bar{p},\bar{x})$. Let $H:P\rightrightarrows P\times X$ be defined by $H(p)=(H_1(p),H_2(p))$, where $H_1(p)=\{p\}$ and $H_2(p)=C(p)$. Suppose that the following conditions
		\begin{enumerate}
			\item[\rm (i)] $\mathbb{R}^+({\rm reg}H-{\rm dom}\,f)$ is a closed subspace of $P\times X$,
			\item[\rm (ii)] $\mathbb{R}^+(P-{\rm dom}\,C)$ is a closed subspace of $P$
		\end{enumerate}
		are satisfied. In addition, assume that $C$ satisfies (ACQ) at $(\bar{p},\bar{x})$ and  $F$ is $K$-dominated by $\mathcal{F}$ near $\bar{p}$. Then, for any $y^*\in Y^*$, we have
\begin{equation*}
	\begin{split}
	\widehat{D}^*(\mathcal{F}\!+\!K)(\bar{p},\bar{y})(y^*) \!=\!\!\bigcup\limits_{(p^*,x^*)\!\in\! D^*(f\!+\!K)((\bar{p},\bar{x}))(y^*)}\bigcup\limits_{\lambda\in \Lambda(\bar{p},\bar{x},x^*)}
	\left\{p^*
\!\!+\!\!\sum\limits_{i=1}^{m+r}\lambda_i\nabla_p g_i(\bar{p},\bar{x})\right\},
\end{split}
\end{equation*}
where
\begin{align*}\Lambda(\bar{p},\bar{x},x^*)\!=\!\bigg\{\lambda=(\lambda_1,...,\lambda_{m+r})\in \mathbb{R}^{m+r}\mid -x^*\!\!=\!\!\sum\limits_{i=1}^{m+r}\lambda_i\nabla_x g_i(\bar{p},\bar{x}),\\
	\lambda_i \ge 0, \, \lambda_i g_i(\bar p, \bar x)=0 \ \mbox{\rm for}\ i=1,...,m
\bigg\}.
\end{align*}
In particular, if $f$ is continuously differentiable at $(\bar{p},\bar{x})$, then for all $y^*\in K^*$, 
\begin{equation*}
	\begin{split}
	\widehat{D}^*(\mathcal{F}\!+\!K)(\bar{p},\bar{y})(y^*)\!=\!\! \bigcup\limits_{\lambda \in \Lambda\big(\bar{p},\bar{x},\nabla_x f(\bar p, \bar x)^* (y^*)\big)}\!\!\bigg\{\!\nabla_p f(\bar{p},\bar{x})^*(y^*)\!\!+\!\!\sum\limits_{i=1}^{m+r}\lambda_i\nabla_p g_i(\bar{p},\bar{x})\!\bigg\}.
		\end{split}
\end{equation*}

\end{theorem}
\begin{proof}The proof is immediate from Theorem \ref{Thm_frontier_map} and Proposition \ref{prop4.1} and so it is omitted.
		\qed
\end{proof}
\subsection{Constraints Described by an Arbitrary (Possibly Infinite) Number of Inequalities}
We consider problem {(${\rm PSVO}_p)$ }with the functional constraints described by an arbitrary (possibly infinite) number of inequalities given as follows
\begin{align}\label{Semi_contraint}C(p)=\{x\in X\mid g_t(p,x)\le 0, \ t\in T\},\end{align}
where {$T$ is an arbitrary (possibly infinite) index set and for each $t \in T$ the function $g_t:P\times X\to \overline{\mathbb{R}}$ are assumed to be convex}. {Constraints of type~\eqref{Semi_contraint} are known as \textit{semi-infinite/infinite inequality constraints}. Semi-infinite optimization programming and its wide applications have attracted much attention from many researchers. We refer the reader to the book of Goberna and \O \ \cite[Chapter 5]{Goberna_Lopez_b}, some papers \cite{C11,Chuong_Huy_Yao} and the references therein for more details and discussions in this topic.}

It should be noted that $C$ is a convex set-valued map since $g_t (t\in T)$ are convex. In addition, $C$ is closed as well. Designate $\mathbb{R}^{|T|}_+$ the set of all the functions $ \lambda: T \to  \mathbb{R}$ taking values $\lambda_t$ positive only at finitely many points of $T$, and equal to zero at the other points. Denote by  $$A(\bar{p},\bar{x}):=\{\lambda\in \mathbb{R}^{|T|}_+\mid \lambda_tg_t(\bar{p},\bar{x})=0, \ t\in T\}$$
the set of \textit{active constraint multipliers}.

Let $(\bar{p},\bar{x})\in {\rm gph}\,C$. We say that $C$ satisfies the\textit{ basic constraint qualification} (BCQ) at $(\bar{p},\bar{x})$  if
$$
N((\bar{p},\bar{x}), {\rm gph}\,C)\subset \bigcup\limits_{\lambda\in A(\bar{p},\bar{x})}\bigg[\sum\limits_{t\in T}\lambda_t\partial g_t(\bar{p},\bar{x})\bigg].
$$
The basic constraint qualification was introduced by Hiriart Urruty~\cite[p.~307]{Hiriart_Urruty_1993} in relation to an ordinary convex programming problem, with equality/inequality constraints.
Various criteria for the validity of this qualification condition can be found in~\cite[Theorems 3.2 and 3.4]{Chuong_Huy_Yao} and \cite[Corollary 2]{Dinh_Mordukhovich_Nghia}.

\medskip

 We are now in a position to establish formulas for computing the Fr\'{e}chet coderivative of the profile of $\mathcal{F}$ in a semi-infinite vector optimization problem.
\begin{theorem}\label{Thm6} Let $f:P\times X\to Y$ be $K$-closed and $K$-convex, $C:P\rightrightarrows X$ be a convex closed set-valued map given by~\eqref{Semi_contraint}, $(\bar{p},\bar{x})\in {\rm gph}\,\mathcal{S}$ and $\bar{y}=f(\bar{p},\bar{x})$. Let $H:P\rightrightarrows P\times X$ be defined by $H(p)=(H_1(p),H_2(p))$, where $H_1(p)=\{p\}$ and $H_2(p)=C(p)$. Suppose that the following conditions
		\begin{enumerate}
			\item[\rm (i)] $\mathbb{R}^+({\rm reg}H-{\rm dom}\,f)$ is a closed subspace of $P\times X$,
			\item[\rm (ii)] $\mathbb{R}^+(P-{\rm dom}\,C)$ is a closed subspace of $P$
		\end{enumerate}
		are satisfied. In addition, assume that $C$ satisfies (BCQ) at $(\bar{p},\bar{x})$ and $F$ is $K$-dominated by $\mathcal{F}$ near $\bar{p}$.
 Then,  for all $ y^*\in Y^*$,
 \begin{equation*}
 	\begin{split}
 	&	\widehat{D}^*(\mathcal{F}\!+\!K)(\bar{p},\bar{y})(y^*) \\
 		&\quad \quad =\bigcup\limits_{(p^*,x^*)\!\in\! D^*(f\!+\!K)((\bar{p},\bar{x}))(y^*)}\bigg\{\!p^*
 		+ u^*\!\mid\! (u^*,\!-\!x^*)\!\in\! \!\bigcup\limits_{\lambda\in A(\bar{p},\bar{x})}\!\sum\limits_{t\in\! T}\lambda_t\partial g_t(\bar{p},\bar{x})\!\bigg\}.
 	\end{split}
 \end{equation*}
In particular, if $f$ is continuously differentiable at $(\bar{p},\bar{x})$, then for  all $ y^*\in K^*$,
 \begin{equation*}
	\begin{split}
		&\widehat{D}^*(\mathcal{F}\!+\!K)(\bar{p},\bar{y})(y^*)\\
		&\quad \quad= \bigg\{\!\nabla_p f(\bar{p},\bar{x})^*(y^*)\!+\!u^*\!\mid\! (u^*,\!-\!\nabla_x f(\bar{p},\bar{x})^*(y^*))\!\in\! \!\bigcup\limits_{\lambda\in A(\bar{p},\bar{x})}\!\sum\limits_{\!t\in\! T}\lambda_t\partial g_t(\bar{p},\bar{x})\!\bigg\}.	\end{split}
	\end{equation*}
\end{theorem}
\begin{proof} Since $C$ satisfies  (BCQ) at $(\bar{p},\bar{x})$, one gets
$$D^*C(\bar{p},\bar{x})(x^*)=\bigg\{{u^*\in P^*}\mid (u^*,-x^*)\in \bigcup\limits_{\lambda\in A(\bar{p},\bar{x})}\sum\limits_{t\in T}\lambda_t\partial g_t(\bar{p},\bar{x})\bigg\}.$$
By invoking Theorem \ref{Thm_frontier_map}, the conclusions are obtained.
	\qed
\end{proof}
\medskip

 In the last part of the paper, we will present an example to illustrate the established results.
\begin{example} Let $T=[0,1]$, $P=\mathbb{R}, X=Y=\mathbb{R}^2$, $K=\mathbb{R}^2_+$, $f:\mathbb{R}\times \mathbb{R}^2\to \mathbb{R}^2$, and $g_t:\mathbb{R}\times \mathbb{R}^2\to \mathbb{R}$ be maps
which are given as follows:
$$f(p,x)=(x_1+2,x_2+3), \forall x=(x_1,x_2)\in\mathbb{R}^2,\forall p\in \mathbb{R},$$
$$g_t(p,x)=-tx_1-(1-t)x_2, , \forall x=(x_1,x_2)\in\mathbb{R}^2,\forall p\in \mathbb{R}.$$
Then, $g_t$ is convex at any $(p, x) \in P \times X$ for all $t\in T$. By some calculation, one has
$$C(p)=\{x\in \mathbb{R}^2\mid x_1\ge 0, x_2\ge 0\},$$
$$F(p)=\{y\in \mathbb{R}^2\mid y_1\ge 2, y_2\ge 3\}.$$
Hence,
$$\mathcal{F}(p)=\{(2,3)\},$$
$$\mathcal{S}(p)=\{x\in X\mid x\in C(p), f(p,x)\in \mathcal{F}(p)\}=\{(0,0)\}.$$
We see that $f$ is $K$-convex, $C$ is convex, the domination property holds for $F$ at all $p \in P$ and $\mathcal{F}$ is convex. For $\bar{p}=0\in P$, $\bar{x}=(0,0)\in \mathcal{S}(\bar{p})$ and $\bar{y}=f(\bar{p},\bar{x})=(2,3)$, {we have}
$$\nabla_p f(\bar{p},\bar{x})=[0\quad 0],
\nabla_x  f(\bar{p},\bar{x})=\left[
\begin{array}{ll}
1&0\\0&1
\end{array}
\right],$$
$$\nabla_p f(\bar{p},\bar{x})^*=[0\quad 0]^T=\left[
\begin{array}{l}
0\\0
\end{array}
\right],
\nabla_x  f(\bar{p},\bar{x})^*=\left[
\begin{array}{ll}
1&0\\0&1
\end{array}
\right]^T=\left[\begin{array}{ll}
1&0\\0&1
\end{array}
\right],$$
$${T}((\bar{p},\bar{x}),{\rm gph}\,C)=\mathbb{R}\times\mathbb{R}^2_+,$$
$$N((\bar{p},\bar{x}),{\rm gph}\,C)=T((\bar{p},\bar{x}),{\rm gph}\,C)^-=\{0\}\times (-\mathbb{R}^2_+),$$
$$D^*C(\bar{p},\bar{x})(x^*)=\left\{
\begin{array}{ll}
0,& \mbox{if}\; x^*\in\mathbb{R}^2_+ \\ \emptyset & \mbox{otherwise,}
\end{array}
\right.$$
$$T(\bar{p},\bar{x})=T=[0,1],A(\bar{p},\bar{x})=\mathbb{R}^{|T|}_+,$$
$$\partial g_t(\bar{p},\bar{x})=\left(0,\left[
\begin{array}{l}
-t\\t-1
\end{array}
\right]\right),\forall t\in T=[0,1].$$
Hence,
$$\bigcup\limits_{\lambda\in A(\bar{p},\bar{x})}\sum\limits_{t\in T}\lambda_t\partial g_t(\bar{p},\bar{x})=\{0\}\times (-\mathbb{R}^2_+),$$
i.e., (BCQ) holds. On one hand
$${\rm reg}H=\bigcup\limits_{q\in P}H(q)=\{(p,x)\in P\times X\mid p=q,x_1\ge 2,x_2\ge 3\},$$
and ${\rm dom}\,f=\mathbb{R}^3.$ So
 $\mathbb{R}^+({\rm reg}H-{\rm dom}\,f)=\mathbb{R}^3$ is a closed subspace of $P\times X$, which in turn implies that~(i) is fulfilled. On the other hand, as ${\rm dom}\,C=P$, it follows that  $\mathbb{R}^+(P-{\rm dom}\,C)=P$ is a closed subspace of $P$. In other words~(ii) is satisfied.

\noindent Therefore, all assumptions of
Theorem~\ref{Thm6} are satisfied. Moreover, for every $y^*=(y^*_1,y^*_2)\in K$,
$$\nabla_p f(\bar{p},\bar{x})^*(y^*)=\langle (0,0),(y^*_1,y^*_2) \rangle =\bigg\langle \left[
\begin{array}{l}
0\\0
\end{array}
\right],\left[
\begin{array}{l}
y^*_1\\y^*_2
\end{array}
\right]\bigg\rangle=0,$$
$$-\nabla_x f(\bar{p},\bar{x})^*(y^*)=-\left[\begin{array}{ll}
1&0\\0&1
\end{array}
\right].\left[
\begin{array}{l}
y^*_1\\y^*_2
\end{array}
\right]=\left[
\begin{array}{l}
-y^*_1\\-y^*_2
\end{array}
\right]\in -\mathbb{R}^2_+.$$
Thus, for all $y^*=(y^*_1,y^*_2)\in K$,
$$(u^*,-\nabla_x f(\bar{p},\bar{x})^*(y^*))\in \bigcup\limits_{\lambda\in A(\bar{p},\bar{x})}\sum\limits_{t\in T}\lambda_t\partial g_t(\bar{p},\bar{x})=\{0\}\times (-\mathbb{R}^2_+)$$
leads that $u^*=0$, and hence,
\begin{align*}D^*\mathcal{F}(\bar{p},\bar{y})(y^*)&\!\!=\!\!\left\{\!\nabla_p f(\bar{p},\bar{x})^*y^*\!+\!u^*\!\mid\! (u^*,-\nabla_x f(\bar{p},\bar{x})^*y^*)\!\in\!\!\! \bigcup\limits_{\lambda\in A(\bar{p},\bar{x})}\sum\limits_{t\in T}\lambda_t\partial g_t(\bar{p},\bar{x})\!\!\right\}\\
&=\{0\}.
	\end{align*}
Therefore, the conclusion of Theorem~\ref{Thm6} is valid.  Furthermore, as $\mathcal{F}(p)~=~\{(2,3)\}$, we can check directly that $D^*\mathcal{F}(\bar{p},\bar{y})(y^*)=\{0\}$.

\end{example}
\section{Concluding Remarks}
\markboth{\centerline{\it Sensitivity Analysis  in Parametric  Convex
		Vector Optimization}}{\centerline{\it D.T.V.~An and L.T.~Tung}} \setcounter{equation}{0}

In this paper, the formulas for computing the Fr\'{e}chet coderivative of the  profile of perturbation, weak perturbation, and proper perturbation maps in parametric convex  vector optimization are studied. Because of the convexity assumptions, the conditions set are fairly simple if compared to those in~\cite{CY13}. In addition, our conditions are stated directly on the data of the problem. It is worth emphasizing that our approach is based on convex analysis tools which is different from the approach in~\cite{CY13}.

\medskip

{\small \noindent{\bf Acknowledgements} A part of this work was completed during a stay of the authors at the Vietnam Institute for Advanced Study in Mathematics (VIASM). The authors would like to thank VIASM for its hospitality and support.
}


\begin{thebibliography}{99}
	\bibitem{An_Cesar_2021} An, D.T.V.,  Guti\'{e}rrez, C.:	Differential stability properties in convex scalar and vector optimization. Set Valued Var. Anal.  \textbf{29}, 893--914 (2021)

\bibitem{An_Yen_2015} An, D.T.V., Yen, N.D.: Differential stability of convex optimization problems under inclusion constraints. Appl. Anal. {\bf 94}, 108--128 (2015)

 \bibitem{Aubin_Ekeland_1984} Aubin, J.P., Ekeland, I.:  Applied Nonlinear Analysis, John Wiley, New York, New York (1984)
\bibitem{Aubin_Frankowska_1990} Aubin, J.P., Frankowska, H.:  Set-Valued Analysis. Birkh${\rm \ddot{a}}$user, Boston (1990)

\bibitem{Bao_Mordukhovich}Bao, T.Q., Mordukhovich, B.S.: Necessary conditions for super minimizers in constrained multiobjective optimization. J. Glob. Optim. 43(4), 533–552 (2009)
\bibitem{Bartl}Bartl D.: A short algebraic proof of the Farkas lemma. SIAM J. Optim. \textbf{19}, 234--239 (2008)

\bibitem{Bao_Mordukhovich}Bao, T.Q., Mordukhovich, B.S.: Necessary conditions for super minimizers in constrained multiobjective optimization. J. Glob. Optim. \textbf{43(4)}, 533–552 (2009)

\bibitem{Bonnans_Shapiro_2000}Bonnans, J.F., Shapiro, A.: Perturbation Analysis of Optimization Problems. Springer, New York (2000)


\bibitem{C11} Chuong,  T.D.: Clarke coderivatives of efficient point multifunctions in parametric vector optimization. Nonlinear Anal. {\bf 74}, 273--285 (2011)

\bibitem{Chuong_Huy_Yao} Chuong,  T.D, Huy,.N.Q., Yao, J.-C.:Subdifferentials of marginal functions in semi-infinite programming. SIAM J. Optim. \textbf{3}, 1462--1477 (2009)

\bibitem{CY09}  Chuong,  T.D.,  Yao, J.-C.: Coderivatives of efficient point multifunctions in parametric vector optimization. Taiwanese J. Math. {\bf 13}, 1671--1693 (2009)

\bibitem{CY10}  Chuong, T.D.,  Yao J.-C.: Generalized Clarke epiderivatives of parametric vector optimization problems. J. Optim. Theory Appl. \textbf{146},  77--94 (2010)
\bibitem{CY13}  Chuong, T.D., Yao J.-C.:  \F \  subdifferentials of efficient point multifunctions in parametric vector optimization.  J. Glob. Optim.  \textbf{57}: 1229--1243 (2013)
\bibitem{C13} Clarke, F.: Functional analysis, Calculus of Variations and Optimal control. Springer, New York (2013)
\bibitem{Dinh_Mordukhovich_Nghia}  Dinh, N.,  Mordukhovich, B.S.,  Nghia, T.T.A.: Qualification and optimality conditions for DC programs with infinite constraints, Acta Math. Vietnam. \textbf{34} 123–153 (2009) 

\bibitem{Fiacco_1983} Fiacco, A.V.: Introduction to Sensitivity and Stability Analysis, Academic Press, New York (1983).
\bibitem{Goberna_Lopez_b}  Goberna, M.A.,  \O, \ M.A.: Linear Semi-Infinite Optimization, Wiley, Chichester (1998)
\bibitem{Hiriart_Urruty_1993} Hiriart Urruty, J.-B., Lemarechal, C.: Convex Analysis and Minimization Algorithms I. Springer-Verlag, Berlin (1993)

\bibitem{HMY08}Huy, N.Q., Mordukhovich, B.S., Yao, J.-C.: Coderivatives of frontier and solution maps in parametric multiobjective optimization, Taiwanese J. Math., \textbf{12(8)}, 2083--2111 (2008)

\bibitem{Khan_Tammer_Zanilescu_2015}   Khan, A.A.,   Tammer, C., Z\u{a}nilescu, C. : Set-Valued Optimization, Springer-Verlag, Berlin (2015)

\bibitem{Kuk_Tanino_Tanaka_96} Kuk, H., Tanino, T.,  Tanaka, M.: Sensitivity analysis in parametrized convex vector optimization. J. Math. Anal. Appl. \textbf{202(2)}, 511-522 (1996)

\bibitem{LH07} Lee, G.M., Huy, N.Q.: On sensitivity analysis in vector optimization. Taiwanese J. Math.  \textbf{11(3)}, 945--958 (2007)

\bibitem{WuLi97} Li, W.: Abadie's constraint qualification, metric regularity, and error bounds for differentiable convex inequalities. SIAM J. Optim. \textbf{7(4)}, 966–978 (1997)


\bibitem{Li_Penot_Xue}Li, S., Penot, J.-P., Xue, X.: Codifferential calculus. Set Valued Var. Anal. \textbf{19(4)}, 505--536 (2011)



\bibitem{Luc_book_89}  Luc, D.T.: Theory of Vector Optimization, Springer-Verlag (1989)

\bibitem{Malanowski_1987}  Malanowski, K.: Stability of Solutions to Convex Problems of Optimization, Springer, Berlin (1987)

\bibitem{Mordukhovich_1980}  Mordukhovich, B.S.: Metric approximations and necessary optimality conditions for general classes of
nonsmooth extremal problmes. Sov. Math. Dokl. 22, 526–530 (1980)

\bibitem{Mordukhovich_2006}    Mordukhovich, B.S.: Variational Analysis and Generalized Differentiation. Volume I: Basic Theory, Volume II: Applications. Springer, Berlin (2006)

\bibitem{Mordukhovich_Nam_2022} Mordukhovich, B.S, Nam, N.M.: Convex Analysis and Beyond: Volume I: Basic Theory. Springer Series in Operations Research and Financial Engineering (2022)

\bibitem{MNR17} Mordukhovich, B.S., Nam, N.M., Rector, R.B., Tran, T.: Variational geometric approach to generalized differential and conjugate calculi in convex analysis. Set-Valued  Var. Anal. \textbf{25(4)}, 731--755 (2017)


\bibitem{Rockafellar_Wets_2004}  Rockafellar R.T., Wets, R.J.-B.: Variational Analysis, Springer, Berlin-Heidelberg (2009)

\bibitem{Sawaragi_Nakayama_Tanino_85}    Sawaragi, Y.,  Nakayama, H.,   Tanino, T.:  Theory of Multiobjective Optimization, Academic Press, New York
(1985)

\bibitem{S93} Shi, D.S.:  Sensitivity analysis in convex vector optimization. J. Optim. Theory Appl. \textbf{77(1}), 145--159 (1993)

\bibitem{Taa_2003} Taa, A.:  Subdifferentials of multifunctions and Lagrange multipliers for multiobjective optimization problems. J. Math. Anal. Appl. \textbf{283}, 398--415 (2003)


\bibitem{Taa_2005}Taa, A.: $\varepsilon$-Subdifferentials of set-valued maps and $\varepsilon$-weak Pareto optimality for multiobjective optimization. Math. Methods Oper. Res. \textbf{62}, 187--209 (2005)

	\bibitem{Tanino_1988a} Tanino, T.: Sensitivity analysis in multiobjective optimization. J. Optim. Theory Appl. \textbf{56}, 479--499 (1988)
	
	\bibitem{Tanino_1988b}  Tanino, T.: Stability and sensitivity analysis in convex vector optimization. SIAM J. Control Optim. \textbf{26(3)}, 521--536 (1988)



\bibitem{Tanino_1990}  Tanino, T.: Stability and sensitivity analysis in multiobjective nonlinear programming. Ann. Oper. Res. \textbf{27}, 97--114 (1990)


\bibitem{Tung2016}Tung, L.T.: On higher-order adjacent derivative of perturbation map in parametric vector optimization. J. Inequal. Appl. Article ID 112 (2016)


\bibitem{Tung2020}   Tung, L.T: On higher-order proto-differentiability of perturbation maps. Positivity \textbf{24,} 441–462 (2020)

\bibitem{Tung_2021} Tung, L.T.:  On higher-order proto-differentiability and higher-order asymptotic proto-differentiability of weak perturbation maps in parametric vector optimization.  Positivity \textbf{25}, 579--604 (2021)

\end{thebibliography}
\end{document}